\documentclass[11pt]{article}
\usepackage{iftex}
\ifPDFTeX
\fi
\ifXeTeX
  \usepackage{fontspec}
\fi
\usepackage{amsmath,amssymb,amsthm,booktabs,graphicx}
\usepackage[a4paper,margin=23mm]{geometry}
\usepackage{xcolor}
\usepackage[colorlinks=true,allcolors=blue!45!black]{hyperref}
\newtheorem{theorem}{Theorem}
\newtheorem{proposition}[theorem]{Proposition}
\newcommand{\email}[1]{\href{mailto:#1}{\nolinkurl{#1}}}
\hypersetup{pdftitle={Nine hyperbolic limit cycles in a three-dimensional competitive Lotka-Volterra system},
 pdfauthor={Haoyu Albert Wang, Yuxiu Shao, Wei P. Dai}}
\title{Nine hyperbolic limit cycles in a three-dimensional\\ competitive Lotka--Volterra system}
\author{Haoyu Albert Wang$^{1,2}$\qquad Yuxiu Shao$^3$\qquad Wei P. Dai$^1$\\[.6em]
  \small\begin{tabular}{c}
    $^1$Research Institute of Intelligent Complex Systems, Fudan University\\
    $^2$Institute of Science and Technology for Brain-Inspired Intelligence, Fudan University\\
    $^3$Universit\'e C\^ote d'Azur, Laboratoire J.\ A.\ Dieudonn\'e, Parc Valrose,\\
    06108 Nice Cedex 2, France\\[3pt]
    \footnotesize\email{haoyuwang18@fudan.edu.cn}\quad
    \email{Yuxiu.SHAO@univ-cotedazur.fr}\quad\email{weidai@fudan.edu.cn}
  \end{tabular}}
\date{}
\begin{document}
\maketitle
\begin{abstract}
Gyllenberg and Yan conjectured that a three-dimensional competitive
Lotka--Volterra system can admit at least five limit cycles. We resolve this
conjecture by constructing an explicit strictly competitive quadratic
Lotka--Volterra system in $\mathbb{R}_{>0}^3$ with at least nine geometrically
distinct hyperbolic limit cycles at a single rational interaction matrix
$A_{\mathrm{IX}}$. Seven small-amplitude cycles surround the positive
equilibrium and coexist with two finite-amplitude cycles, yielding five
attracting and four saddle cycles with pairwise disjoint certified minimal
periods. Locally, evaluating the three-dimensional Poincar\'e--Lyapunov scalar
recursion in arbitrary-precision ball arithmetic with forward automatic
differentiation and Krawczyk contraction, we bypass symbolic polynomial
elimination to certify weak foci of exact orders $m=5,6,7,8$ in Zeeman's
class~27 with sharp local cyclicity $m$ (and $m-1$ on the Hopf hypersurface),
together with explicit rational matrices $A_5,A_6,A_7,A_8$ possessing at least
$5,6,7,8$ distinct hyperbolic limit cycles. The same method establishes at
least $5$, $5$, and $6$ coexisting hyperbolic limit cycles in Zeeman's
classes~26, 28, and~29; in class~27 on the Hopf hypersurface, we further show
that order-eight focal vanishing fixes neither boundary heteroclinic neutrality
nor its stability sign, and certify an order-seven weak focus with independent
lower-order focal coordinates inside the neutral-boundary locus. We thank
GPT-6 Astra for helping us compare earlier constructions, identify the main
obstacle to extending them, and develop the present approach.
\end{abstract}

\section{Introduction and main results}\label{sec:intro}
Consider the three-species strictly competitive quadratic Lotka--Volterra system
$\dot z_i=z_i\bigl(r_i-\sum_{j=1}^3\tilde a_{ij}z_j\bigr)$ on $\mathbb{R}_{>0}^3$
with intrinsic growth rates $r_i>0$ and interaction coefficients $\tilde a_{ij}>0$.
Whenever a positive equilibrium $z^*\in\mathbb{R}_{>0}^3$ exists, rescaling
$x_i=z_i/z_i^*$ places the interior equilibrium at $\mathbf{1}=(1,1,1)^T$ and
brings the system to the normalized form
\begin{equation}\label{eq:lv}
 \dot x_i=x_i\sum_{j=1}^3 a_{ij}(1-x_j)
         =-x_i\sum_{j=1}^3 a_{ij}(x_j-1),\qquad
 A=(a_{ij})_{1\le i,j\le 3}=\begin{pmatrix}a&b&c\\d&e&f\\g&h&i\end{pmatrix},\qquad a_{ij}>0,
\end{equation}
with growth-rate vector $r=A\mathbf{1}$. By Hirsch's monotone flow theory
\cite{Hirsch1988}, every nontrivial trajectory in $\mathbb{R}_{>0}^3$ is
asymptotic to an invariant two-dimensional topological carrying simplex $\Sigma$
homeomorphic to the standard probability simplex via radial projection, and
every positive periodic orbit lies in $\Sigma$. Zeeman \cite{Zeeman1993}
partitioned generic three-species competitive Lotka--Volterra systems into $33$
stable nullcline equivalence classes, ruling out periodic orbits in
classes~1--25, establishing Hopf bifurcations in classes~26--31, and showing
that a boundary heteroclinic three-cycle on $\partial\Sigma$ occurs uniquely in
class~27 \cite{MayLeonard1975,Zeeman1993}; van den Driessche and Zeeman
\cite{VanDenDriesscheZeeman1998} subsequently excluded periodic orbits in
classes~32--33, confirming that isolated periodic orbits can occur only in
classes~26--31.

Determining how many positive limit cycles can coexist in \eqref{eq:lv} is a
central problem in nonlinear dynamics and theoretical ecology. Hofbauer and So
\cite{HofbauerSo1994} constructed two limit cycles in class~27 and asked whether
two could be the global maximum; this was disproved by three-cycle examples in
class~29 \cite{GyllenbergYanWang2006} and classes~30--31 \cite{LianLuLuo2008},
followed by four-cycle constructions: Gyllenberg and Yan \cite{GyllenbergYan2009}
proposed three small Hopf cycles plus one boundary-trapped cycle in class~27,
after which Yu, Han and Xiao \cite{YuHanXiao2016} identified a focal-value sign
error in the primary example of \cite{GyllenbergYan2009}, validated the
four-cycle conclusion via a second parameter set in \cite{GyllenbergYan2009},
and constructed four purely local small cycles in classes~27 and~26; four cycles
were also obtained in class~29 \cite{WangHuangWu2011} and class~28
\cite{HuLuLuo2026} (see also \cite{XiaoLi2000,LuLuo2002,JiangEtAl2021,LiJiang2024}).
In 2009, Gyllenberg and Yan \cite{GyllenbergYan2009} conjectured that
\eqref{eq:lv} can admit at least five limit cycles (recorded as
Conjecture~1.2(i) in \cite{YuHanXiao2016}); Tian and Yu \cite{TianYu2013}
further discussed a possible route to eight cycles via seven small cycles and one
outer cycle while noting its computational intractability. Higher counts were
achieved only by relaxing the model class---such as seven small cycles in a
sparse non-Lotka--Volterra quadratic field \cite{TianYu2014}, five or six cycles
in noncompetitive Lotka--Volterra systems \cite{WangEtAl2021,HanLlibreTian2020},
and eight small cycles in a cubic noncompetitive Kolmogorov system
\cite{LuEtAl2025}---leaving four limit cycles as the largest certified count
within the strictly competitive quadratic class \eqref{eq:lv}.

In this paper we construct explicit systems \eqref{eq:lv} with five, six, seven,
eight, and nine coexisting hyperbolic limit cycles. For any positive $T$-periodic
orbit of \eqref{eq:lv} with $\det A\ne 0$, Volterra's time-average identity
$\frac{1}{T}\int_0^T x(t)\,dt=\mathbf{1}$
\cite{HofbauerSigmund1998,CostePeyraudCoullet1979} implies that the vector field
$F_{\mathrm{LV}}(x)=x\odot A(\mathbf{1}-x)$ has negative average divergence
$-\operatorname{tr}A<0$, so the two nontrivial (transverse) Floquet multipliers
$\mu_+=\mu_{\mathrm{slow}}$ and $\mu_-=\mu_{\mathrm{strong}}$ satisfy
$\mu_+\mu_-=\exp(-T\operatorname{tr}A)\in(0,1)$. Moreover, because the
time-reversed variational equation has strictly positive off-diagonal entries
$x_i a_{ij}>0$ ($i\ne j$), Kamke's theorem \cite{Hirsch1988} and
Perron--Frobenius theory ensure that the backward monodromy matrix
$D\varphi_T(x(0))^{-1}$ is entrywise positive with simple dominant eigenvalue
$\mu_-^{-1}>1$ and positive remaining eigenvalues $1$ and $\mu_+^{-1}>0$. Thus
$\mu_+,\mu_-$ are always real and strictly positive, and every hyperbolic
positive limit cycle ($\mu_+\ne 1$) is either attracting ($0<\mu_-<\mu_+<1$) or
of saddle type ($0<\mu_-<1<\mu_+$).

\begin{theorem}[Nine coexisting limit cycles at a single matrix]\label{thm:nine}
Let $A_{\mathrm{IX}}\in\mathbb{Q}_{>0}^{3\times 3}$ be the exact finite-decimal
matrix given in Appendix~\ref{sec:exact} (belonging to Zeeman's class~27). Then
system \eqref{eq:lv} with $A=A_{\mathrm{IX}}$ possesses at least nine
geometrically distinct positive hyperbolic limit cycles in $\mathbb{R}_{>0}^3$,
consisting of seven small-amplitude cycles near $\mathbf{1}$ and two
finite-amplitude cycles. Ordered by strictly increasing minimal positive period
$T_1<T_2<\cdots<T_9$, with $\min_{1\le k\le 8}(T_{k+1}-T_k)>2.8815\times10^{-6}$,
the nine cycles strictly alternate between attracting and saddle stability,
beginning with an attracting cycle $T_1$ and ending with an attracting outer
cycle $T_9$ (five attracting and four saddle cycles in total). In particular,
the supremum $H_{\mathrm{LV3}}$ of coexisting positive limit cycles over all
systems \eqref{eq:lv} satisfies $H_{\mathrm{LV3}}\ge 9$.
\end{theorem}

\begin{figure}[tbp]
\centering
\includegraphics[width=\linewidth]{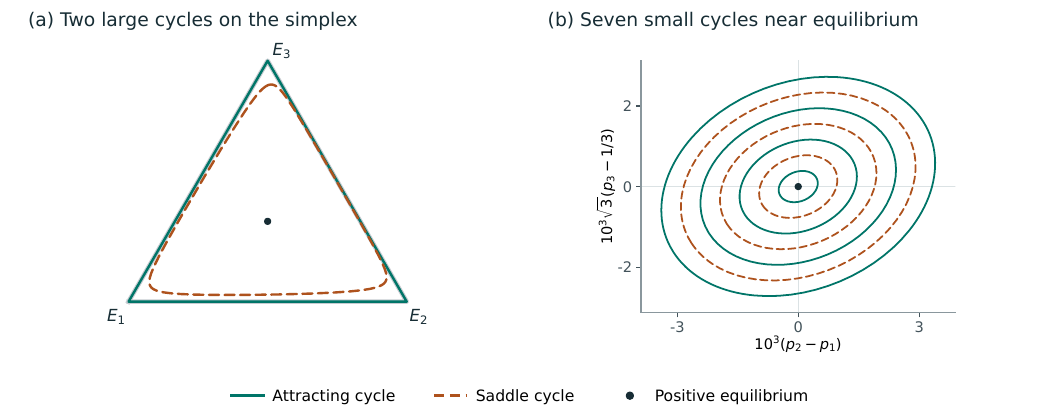}
\caption{Trajectories of the nine certified limit cycles of the fixed matrix
$A_{\mathrm{IX}}$, projected onto the standard simplex in normalized abundance
coordinates $p_i=x_i/(x_1+x_2+x_3)$. (a)~The intermediate saddle cycle (Orbit~8)
and near-boundary attracting cycle (Orbit~9) on the simplex with axial vertices
$E_1,E_2,E_3$. (b)~Magnification by $10^3$ near the positive equilibrium
$\mathbf{1}$, resolving the seven small-amplitude cycles (Orbits~1--7).}
\label{fig:nine-cycles}
\end{figure}

\begin{table}[tbp]
\centering\small
\caption{The nine certified positive hyperbolic limit cycles of the single
rational matrix $A_{\mathrm{IX}}$, ordered by minimal period $T_k$. Displayed
values are rounded midpoints of rigorous interval enclosures: enclosure radii for
$T_k$ are at most $5.94\times 10^{-62}$ (Orbits~1--7) and $3.80\times 10^{-25}$
(Orbits~8--9), and radii for $\delta_{\mathrm{F}}=\mu_{\mathrm{slow}}-1$ are at
most $6.30\times 10^{-55}$ (Orbits~1--7) and $4.15\times 10^{-16}$ (Orbits~8--9).}
\label{tab:lv-nine-certified}
\begin{tabular}{@{}clrrll@{}}
\toprule
Orbit $k$ & Regime & Minimal period $T_k$ & $\delta_{\mathrm{F}}=\mu_{\mathrm{slow}}-1$ & $\mu_{\mathrm{strong}}$ & Stability \\
\midrule
1 & small & 9.693858540 & $-6.953\times 10^{-49}$ & $1.062\times 10^{-23}$ & attracting \\
2 & small & 9.693861422 & $+1.043\times 10^{-48}$ & $1.062\times 10^{-23}$ & saddle \\
3 & small & 9.693866224 & $-2.086\times 10^{-48}$ & $1.062\times 10^{-23}$ & attracting \\
4 & small & 9.693872948 & $+5.736\times 10^{-48}$ & $1.062\times 10^{-23}$ & saddle \\
5 & small & 9.693881593 & $-2.295\times 10^{-47}$ & $1.062\times 10^{-23}$ & attracting \\
6 & small & 9.693892158 & $+1.491\times 10^{-46}$ & $1.062\times 10^{-23}$ & saddle \\
7 & small & 9.693904646 & $-2.089\times 10^{-45}$ & $1.062\times 10^{-23}$ & attracting \\
8 & finite-amplitude & 16.575702060 & $+4.454\times 10^{-10}$ & $5.207\times 10^{-40}$ & saddle \\
9 & finite-amplitude & 226.094115684 & $-4.962\times 10^{-10}$ & $1.480\times 10^{-536}$ & attracting \\
\bottomrule
\end{tabular}
\end{table}

Figure~\ref{fig:nine-cycles} and Table~\ref{tab:lv-nine-certified} summarize the
geometry, minimal periods, and transverse Floquet multipliers of the nine
certified orbits of $A_{\mathrm{IX}}$. We also establish a complete sequence of
purely local weak-focus and unfolding results from order five to order eight in
class~27 and across classes~26, 28, and~29.

\begin{figure}[tbp]
 \centering\includegraphics[width=.84\linewidth]{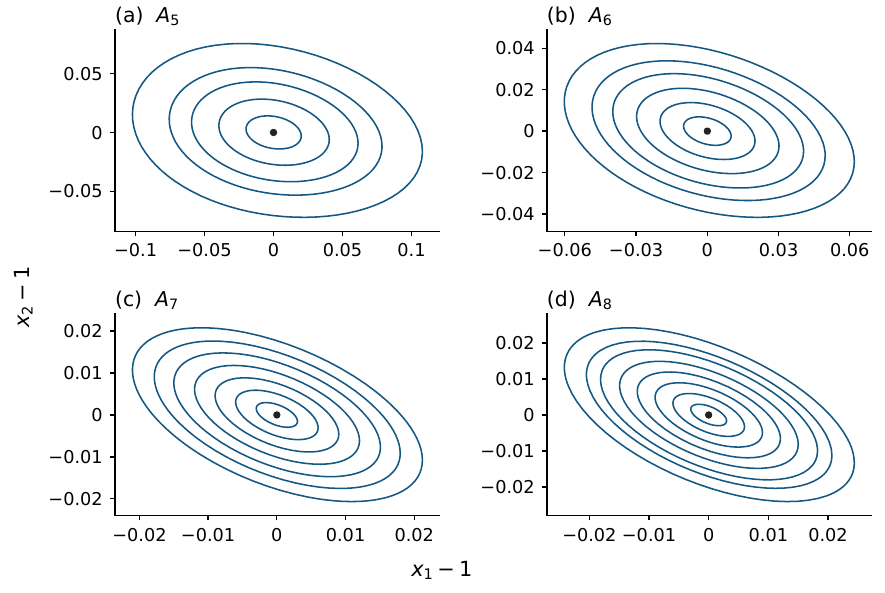}
 \caption{Projections onto $(x_1-1,x_2-1)$ of the certified local limit cycles
 for the fixed rational matrices $A_5,A_6,A_7,A_8$ ($m=5,6,7,8$ cycles). Central
 dots mark $\mathbf{1}$; panel scales differ.}
 \label{fig:orbits}
\end{figure}

\begin{theorem}[Progression from five to eight local cycles and across Zeeman classes]\label{thm:local-and-classes}
\begin{enumerate}
\item[\textup{(a)}] For each $m\in\{5,6,7,8\}$, there exists a strictly positive
critical matrix $A_m^*$ in Zeeman's class~27 on the Hopf hypersurface
$\mathcal{H}$ at which $\mathbf{1}$ is a weak focus of exact order $m$ whose
first $m-1$ focal quantities have full Jacobian rank $m-1$ along $\mathcal{H}$.
The local cyclicity of $\mathbf{1}$ at $A_m^*$ is exactly $m$ in the full family
\eqref{eq:lv} and exactly $m-1$ on $\mathcal{H}$. Moreover, for each
$m\in\{5,6,7,8\}$, there is an explicit rational matrix
$A_m\in\mathbb{Q}_{>0}^{3\times 3}$ (see Appendix~\ref{sec:exact}) defining a
system \eqref{eq:lv} in class~27 with at least $m$ distinct positive hyperbolic
limit cycles: ordered from innermost to outermost, the $m$ certified
small-amplitude cycles near $\mathbf{1}$ alternate between attracting and saddle
stability and end with an attracting cycle (Figure~\ref{fig:orbits}).
\item[\textup{(b)}] In Zeeman's classes~26, 28, and~29, there exist strictly
positive critical matrices on $\mathcal{H}$ that are certified weak foci of
exact orders $5$, $5$, and $6$, respectively, with full-rank $(m-1)$-parameter
unfoldings along $\mathcal{H}$ and transverse Hopf crossing. Consequently, each
of classes~26 and~28 admits systems \eqref{eq:lv} with at least five coexisting
positive hyperbolic limit cycles near $\mathbf{1}$, and class~29 admits systems
with at least six.
\end{enumerate}
\end{theorem}

\section{Local weak foci, parameter selection, and interval validation}\label{sec:method}

Earlier constructions in \eqref{eq:lv} computed symbolic focal quantities on
low-dimensional parameter slices, either combining up to three Hopf cycles with
one boundary-trapped cycle
\cite{HofbauerSo1994,GyllenbergYan2009,WangHuangWu2011,LiJiang2024,HuLuLuo2026}
or solving polynomial elimination systems for four small cycles
\cite{YuHanXiao2016}. Because strict competition ($a_{ij}>0$) prevents
sparsifying the vector field, symbolic focal expressions swell rapidly by fourth
order \cite{TianYu2013,TianYu2014} and two- or three-parameter slices cannot
support full-rank unfoldings of order $m\ge 5$; moreover, qualitative Bautin
unfolding does not certify concrete matrices when slow Floquet offsets reach
$10^{-42}$--$10^{-49}$ (Tables~\ref{tab:lv-nine-certified} and~\ref{tab:fixed},
where even $30$-digit rounding of $A_{\mathrm{IX}}$ flips the sign of
$H_{\mathrm{RH}}$ from $-3.791\times 10^{-48}$ to $+3.664\times 10^{-30}$), and
annular Poincar\'e--Bendixson trapping yields no conclusion when
$L_m P_{\mathrm{het}}<0$. We bypass these obstructions in three steps: direct
ball-arithmetic evaluation of the focal jet and its Jacobian on the full Hopf
hypersurface, multiprecision parameter selection for explicit rational matrices
(including the $7+2$ configuration $A_{\mathrm{IX}}$), and direct interval
validation of all periodic orbits in the original differential equations
\eqref{eq:lv}.

\subsection{Weak foci and sharp local cyclicity on the Hopf hypersurface}\label{sec:focus}
To avoid symbolic polynomial expansion on low-dimensional slices, we evaluate the
focal jet directly in ball arithmetic on a global rational chart of the Hopf
hypersurface $\mathcal{H}$. For $u=x-\mathbf{1}$, system \eqref{eq:lv} becomes
$\dot u=Ju+Q(u)$ with $J=-A$ and $Q(u)=-\operatorname{diag}(u)Au$. The
characteristic polynomial of $J$ is
$\lambda^3+(\operatorname{tr}A)\lambda^2+S\lambda+D$, with
$\operatorname{tr}A=a+e+i$, $S=ae+ai+ei-bd-cg-fh$, and $D=\det A$. On the Hopf
hypersurface $\mathcal{H}=\{A>0:H_{\mathrm{RH}}:=S\operatorname{tr}A-D=0,\ S>0\}$,
$J$ has eigenvalues $\pm\mathrm{i}\omega$ ($\omega=\sqrt{S}>0$) and
$-\operatorname{tr}A<0$. Since $\partial_c(D-S\operatorname{tr}A)=dh+(a+i)g>0$,
the Hopf equation $H_{\mathrm{RH}}=0$ globally determines $c$ as a rational
function of the remaining eight entries, and the real part $\alpha$ of the
oscillatory eigenvalue pair $\alpha\pm\mathrm{i}\omega$ of $J$ satisfies the
strict Hopf transversality condition
\begin{equation}\label{foc-cross}
 \partial_c\alpha
 =\frac{dh+(a+i)g}{2\bigl(S+(\operatorname{tr}A)^2\bigr)}>0.
\end{equation}

To fix the eigenbasis used in the computation, take the cross product of the
first two rows of $A+\mathrm{i}\omega I$ for $q$ and of
$A-(\operatorname{tr}A)I$ for $v_{\rm st}$, divide each vector by its first
component (unit first entry), and let $\bar q$ be the complex conjugate of $q$
(this normalization uniquely fixes the numerical scale of $L_j(A)$ in
Table~\ref{tab:focus}). On every certified parameter box, these normalizing
first components, the chart denominator $dh+(a+i)g$, and
$\det[q\ \bar q\ v_{\rm st}]$ are enclosed strictly away from zero, so
$u=zq+\bar z\,\bar q+wv_{\rm st}$ defines smooth coordinates
$(z,\bar z,w)^T=[q\ \bar q\ v_{\rm st}]^{-1}u$. Among the three equivalent
classical formulations for three-dimensional Hopf focal
quantities---center-manifold power-series reduction \cite{Carr1981,Kuznetsov2004},
direct normal-form recursion \cite{TianYu2013,YuHanXiao2016}, and the
three-dimensional Poincar\'e--Lyapunov formal first-integral recursion
\cite{EdneralEtAl2012,RomanovskiShafer2014,GarciaMazaShafer2018}---we adopt the
third. Let $\mathcal{L}$ denote the linear Lie derivative along $\dot u=Ju$,
which multiplies each monomial $z^j\bar z^k w^\ell$ of total degree
$n=j+k+\ell$ by $\mathrm{i}\omega(j-k)-\ell\operatorname{tr}A$; because
$\omega>0$ and $\operatorname{tr}A>0$, this divisor vanishes if and only if
$j=k$ and $\ell=0$. Let $\Pi_n$ project homogeneous degree-$n$ polynomials
onto the resonant subspace
\[
 \operatorname{im}\Pi_n=\begin{cases}
  \operatorname{span}\{(z\bar z)^{n/2}\}, & n\text{ even},\\
  \{0\}, & n\text{ odd},
 \end{cases}
\]
annihilating all nonresonant monomials (which span $\ker\Pi_n$). Starting from
$V_2(z,\bar z,w)=z\bar z$, we solve for $3\le n\le 2m+2$
\begin{equation}\label{foc-rec}
 R_n=DV_{n-1}\,Q,\qquad
 \mathcal{L}V_n=-(I-\Pi_n)R_n,\qquad \Pi_nV_n=0,
\end{equation}
where $\mathcal{L}|_{\ker\Pi_n}$ is inverted termwise by dividing by the
nonzero divisors $\mathrm{i}\omega(j-k)-\ell\operatorname{tr}A$, and
$\Pi_{2j+2}R_{2j+2}=L_j(A)|z|^{2j+2}$ defines the $j$th scalar focal quantity
$L_j(A)$, giving $\dot V=\sum_{j=1}^m L_j(A)|z|^{2j+2}+O(\|u\|^{2m+3})$ for
$V=\sum_{n=2}^{2m+2}V_n$ (terms involving $w$ are retained at every degree
since their products with $Q$ can contribute to $\operatorname{im}\Pi_{2j+2}$
at higher degrees).

Instead of expanding \eqref{foc-rec} symbolically, we fix $a=1$ and $8-m$
additional entries, vary $m-1$ free entries $\xi\in\mathbb{R}^{m-1}$ on
$\mathcal{H}$, and evaluate $F(\xi)=(L_1(\xi),\dots,L_{m-1}(\xi))$ and its
Jacobian $DF(X)$ over a narrow rational box
$X=\{\xi:\|\xi-\xi_0\|_\infty\le\Delta\}$ directly in outward-rounded
arbitrary-precision ball arithmetic \cite{Johansson2017} via forward automatic
differentiation \cite{GriewankWalther2008,Tucker2011}. With
$R\approx DF(\xi_0)^{-1}$, we verify the Krawczyk contraction inequalities
\cite{Krawczyk1969,Rump2010,Tucker2011}
\[
 q:=\max_{1\le k\le m-1}q_k<1,\qquad
 q_k:=\sup_{\xi\in X}\sum_{\ell=1}^{m-1}|(I-RDF(\xi))_{k\ell}|,\qquad
 \max_{1\le k\le m-1}h_k<\Delta,
\]
where $h_k:=|(RF(\xi_0))_k|+\Delta q_k$, together with $0\notin L_m(X)$,
$A(X)>0$, $S(X)>0$, and nonvanishing of the chart and eigenbasis denominators
across $X$, certifying a unique zero $\xi_*\in\operatorname{int}(X)$ with
nonsingular Jacobian $DF(\xi_*)$ (Table~\ref{tab:focus}).

\paragraph{Sharp local cyclicity via center-manifold reduction and Rolle's theorem.}
Fix any certified order-$m$ critical matrix $A_*\in\mathcal{H}$ from
Table~\ref{tab:focus} and write $\lambda_*:=L_m(A_*)\ne 0$. To relate the
certified focal rank to exact local cyclicity, we combine parameter-dependent
center-manifold reduction \cite{Carr1981,Vanderbauwhede1989,Kuznetsov2004} with
the Bautin--Rolle displacement argument
\cite{Bautin1954,Roussarie1998,Christopher2005,GarciaMazaShafer2018}. Adjoining
$\dot A=0$ to $\dot u=Ju+Q(u)$ near $(0,A_*)$ yields a local $C^{2m+4}$ center
manifold whose constant-$A$ slices $W_A^c=\{w=h_A(z,\bar z)\}$
($h_A(z,\bar z)=O(|z|^2)$) are two-dimensional locally invariant surfaces
containing every periodic orbit of \eqref{eq:lv} in a sufficiently small
neighborhood of $\mathbf{1}$ \cite{Vanderbauwhede1989}. This local argument
uses the smoothness of $W_A^c$ rather than that of Hirsch's global carrying
simplex $\Sigma$, which need only be $C^1$ or Lipschitz in general.

Continuing the eigenbasis smoothly to matrices $A$ near $A_*$, the reduced
vector field on $W_A^c$ has the form $\dot z=(\alpha+\mathrm{i}\omega)z+O(|z|^2)$
with $\omega>0$. Restricting the polynomial $V$ from \eqref{foc-rec} to $W_A^c$
for $A\in\mathcal{H}$ near $A_*$ and setting $s=|z|^2$ gives, uniformly in
$(z,A)$,
\[
 V|_{W_A^c}=s+O(s^{3/2}),\qquad
 \dot V|_{W_A^c}=\sum_{j=1}^m L_j(A)s^{j+1}+O(s^{m+3/2}).
\]
Since $\partial_\rho(V|_{W_A^c})=2\rho+O(\rho^2)>0$ for small $\rho=|z|>0$, the
small positive level sets of $V|_{W_A^c}$ are nested closed ovals transverse to
rays from the origin.

To obtain the local upper bounds, take the real $z$-axis as a section on
$W_A^c$, write $P_A(\sigma)$ for the signed coordinate after one full
revolution $\Delta\theta=2\pi$ starting from $z=\sigma\in\mathbb{R}$, and set
$H_A(\sigma):=P_A(\sigma)-\sigma$ with $P_A(0)=0$. Writing
$z=\varrho e^{\mathrm{i}\theta}$ with signed $\varrho\in\mathbb{R}$, Hadamard's
integral lemma divides the $C^{2m+4}$ reduced field by $\varrho$ with the loss
of one derivative, and since $\dot\theta=\omega+O(|\varrho|)>0$, integrating
$d\varrho/d\theta$ from $0$ to $2\pi$ defines $H_A(\sigma)$ as a $C^{2m+3}$
function of $(\sigma,A)$ near $(0,A_*)$ \cite{Roussarie1998}. At $A_*$, one
revolution takes time $T_*+O(|\sigma|)$ with $T_*=2\pi/\omega_*$, giving
$z(t)=\sigma e^{\mathrm{i}\omega_* t}+O(\sigma^2)$ and
$H_{A_*}(\sigma)=O(\sigma^2)$. Letting $v(\sigma):=V|_{W_{A_*}^c}(\sigma)$ on
the section ($v'(\sigma)=2\sigma+O(\sigma^2)$), integrating $\dot V$ over one
turn and applying the integral mean value identity gives
\begin{align*}
 T_*\lambda_*\sigma^{2m+2}+O(|\sigma|^{2m+3})
 &=v\bigl(\sigma+H_{A_*}(\sigma)\bigr)-v(\sigma)\\
 &=H_{A_*}(\sigma)\int_0^1 v'\bigl(\sigma+\tau H_{A_*}(\sigma)\bigr)\,d\tau
 =\bigl(2\sigma+O(\sigma^2)\bigr)H_{A_*}(\sigma),
\end{align*}
whence
\begin{equation}\label{foc-scalar-return}
 H_{A_*}(\sigma)=\frac{\pi\lambda_*}{\omega_*}\sigma^{2m+1}+O(|\sigma|^{2m+2}).
\end{equation}
Thus $\partial_\sigma^{2m+1}H_A(\sigma)$ remains strictly nonzero on a uniform
interval $[-\sigma_0,\sigma_0]$ for all $A$ sufficiently close to $A_*$. By
repeated application of Rolle's theorem, $H_A$ has at most $2m+1$ zeros on
$[-\sigma_0,\sigma_0]$ counted with multiplicity. Because $P_A$ is
strictly increasing on the two-sided section (so every periodic point of
$P_A$ is fixed), every nontrivial periodic orbit in the neighborhood
meets both half-sections $\sigma>0$ and $\sigma<0$, contributing two nonzero
roots of $H_A$ in addition to $\sigma=0$. Moreover, for $A\in\mathcal{H}$ one
has $\dot V|_{W_A^c}=O(|z|^4)$, so the same mean-value division yields
$H_A(\sigma)=O(|\sigma|^3)$ and $\sigma=0$ has multiplicity at least $3$. Hence
any $N$ coexisting small limit cycles near $\mathbf{1}$ satisfy
$1+2N\le 2m+1$ ($N\le m$) in \eqref{eq:lv} and $3+2N\le 2m+1$ ($N\le m-1$) on
$\mathcal{H}$.

To attain both bounds simultaneously \cite{Christopher2005}, choose $\delta>0$
small and set $Q_\delta(s):=\lambda_* s^2\prod_{k=1}^{m-1}(s-k\delta)$. By the
nonsingular $(m-1)\times(m-1)$ Jacobian certified in Table~\ref{tab:focus} and
the inverse function theorem on the Hopf slice, there exists a unique matrix
$A_\delta=A_*+O(\delta)\in\mathcal{H}$ matching $L_j(A_\delta)$ to the
coefficient of $s^{j+1}$ in $Q_\delta(s)$ for $1\le j\le m-1$, while
$L_m(A_\delta)=\lambda_*+O(\delta)$. On the $m$ nested ovals
$\{V|_{W_{A_\delta}^c}=\delta t\}$ for $t=\tfrac12,\tfrac32,\dots,m-\tfrac12$,
one has $s=\delta t+O(\delta^{3/2})$ and $L_j(A_\delta)=O(\delta^{m-j})$,
whence
\[
 \delta^{-(m+1)}\dot V|_{W_{A_\delta}^c}
 =\lambda_* t^2\prod_{k=1}^{m-1}(t-k)+O(\sqrt{\delta})
\]
uniformly in angle. For small $\delta>0$, $\dot V|_{W_{A_\delta}^c}$ therefore
has strict alternating signs on the $m$ ovals, preventing trajectories from
re-crossing any test oval during one revolution and forcing $H_{A_\delta}$ to
alternate in sign at the $m$ corresponding positive section points.
Intermediate-value zeros give $m-1$ positive fixed points of
$P_{A_\delta}$ on $\mathcal{H}$. Fixing $\delta$ and perturbing
$c$ slightly via \eqref{foc-cross} so that $\alpha$ has the sign opposite
to $H_{A_\delta}$ at the innermost test oval preserves the $m$ oval return
signs while giving $P_A'(0)=\exp(2\pi\alpha/\omega)\ne 1$, creating
an $m$th positive fixed point near $\sigma=0$. Because $1+2m=2m+1$ (and
$3+2(m-1)=2m+1$ on $\mathcal{H}$) exhausts the Rolle multiplicity bound, every
nonzero root of $H_A$ is simple ($P_A'(\sigma)\ne 1$), and the second
transverse Floquet multiplier approaches $\exp(-T_*\operatorname{tr}A_*)<1$.
Thus all $m$ (resp.\ $m-1$) bifurcating cycles are hyperbolic, and entrywise
positivity of $A$ together with the strict invasion inequalities defining the
Zeeman class persists for all sufficiently small perturbations, proving the
local assertions of Theorem~\ref{thm:local-and-classes}(a)--(b).

\begin{table}[tbp]
\centering\footnotesize
\setlength{\tabcolsep}{3.5pt}
\caption{Certified weak-focus parameter boxes on the Hopf hypersurface
$\mathcal{H}$ ($a=1$, $c$ solved from $D=S\operatorname{tr}A$). Here $B_0$
strictly bounds both $q$ and $\max_k h_k/\Delta$ (for $A_{7,\mathrm{IX}}^*$,
$q<1.37\times 10^{-81}$ and $\max_k h_k/\Delta<4.17\times 10^{-20}$), and the
last column encloses $L_m$ over $X$.}
\label{tab:focus}
\begin{tabular}{@{}lclccc@{}}
\toprule
Critical point (class) & Order $m$ & Free entries & $\Delta$ & $B_0$ & $L_m(X)$ \\
\midrule
$A_5^*$ (class~27) & 5 & $e,i,g,b$ & $10^{-30}$ & $3.19\times 10^{-17}$ & $[-2.06392,-2.06391]\times 10^{-6}$ \\
$A_6^*$ (class~27) & 6 & $e,g,d,f,i$ & $10^{-30}$ & $1.36\times 10^{-15}$ & $[-1.83957,-1.83956]\times 10^{-6}$ \\
$A_7^*$ (class~27) & 7 & $e,b,d,g,f,i$ & $10^{-32}$ & $1.62\times 10^{-14}$ & $[-3.79062,-3.79061]\times 10^{-8}$ \\
$A_8^*$ (class~27) & 8 & $e,f,d,g,b,i,h$ & $10^{-32}$ & $7.64\times 10^{-13}$ & $[-2.61971,-2.61970]\times 10^{-7}$ \\
$A_{7,\mathrm{IX}}^*$ (27, 9-cycle) & 7 & $f,b,e,d,i,h$ & $10^{-100}$ & $4.17\times 10^{-20}$ & $[-8.09643,-8.09642]\times 10^{-11}$ \\
Class~26 focus & 5 & $d,b,g,i$ & $10^{-42}$ & $5.76\times 10^{-28}$ & $[-2.60389,-2.60388]\times 10^{-7}$ \\
Class~28 focus & 5 & $e,d,i,h$ & $10^{-42}$ & $2.14\times 10^{-26}$ & $[7.23433,7.23434]\times 10^{-7}$ \\
Class~29 focus & 6 & $e,g,i,d,f$ & $10^{-42}$ & $9.64\times 10^{-27}$ & $[-6.94427,-6.94426]\times 10^{-7}$ \\
\bottomrule
\end{tabular}
\end{table}

\subsection{Choosing the fixed matrices \texorpdfstring{$A_5,\dots,A_8$}{A5,...,A8} and the nine-cycle matrix \texorpdfstring{$A_{\mathrm{IX}}$}{A\_IX}}
To select a fixed matrix off the Hopf hypersurface $\mathcal{H}$, write the
eigenvalues of $J=-A$ as $\alpha\pm\mathrm{i}\omega$ and
$\lambda_{\rm st}=-\operatorname{tr}A-2\alpha$, and set $s=|z|^2$ in the
continued linear eigenbasis. Retaining the same projection $\Pi_n$ onto powers of
$s$ and the same normalization $\Pi_n V_n=0$ in \eqref{foc-rec}, the divisors
$(j+k)\alpha+\mathrm{i}(j-k)\omega+\ell\lambda_{\rm st}$ associated with
monomials $z^j\bar z^k w^\ell$ in $\ker\Pi_n$ ($j\ne k$ or $\ell\ge 1$) remain
nonzero near $A_*$, so $\mathcal{L}_A|_{\ker\Pi_n}$ remains invertible. Hence
the coefficients $L_j(A)$ extend smoothly off $\mathcal{H}$, giving the
off-Hopf scalar drift
\begin{equation}\label{eq:finite-scalar}
 \dot V=sD_A(s)+O(\|u\|^{2m+3}),\qquad
 D_A(s)=\sum_{k=0}^m L_k(A)s^k,\qquad L_0(A)=2\alpha.
\end{equation}
For a scale $\varepsilon>0$, we prescribe target radii $j\varepsilon$
($1\le j\le m$) and match
\begin{equation}\label{eq:finite-target}
 \prod_{j=1}^m(s-j^2\varepsilon^2)=\sum_{k=0}^m p_k(\varepsilon)s^k,\qquad
 L_k(A)=L_m(A)\,p_k(\varepsilon),\quad 0\le k<m,
\end{equation}
solving for the free entries in Table~\ref{tab:focus} together with $\alpha$,
while recovering $c$ explicitly from the characteristic polynomial: with
$S_0=ae+ai+ei-bd-fh$ and $D_0=aei-afh-bdi+bfg$,
\[
 c=\frac{(\operatorname{tr}A+2\alpha)S_0+2\alpha(\operatorname{tr}A+2\alpha)^2-D_0}{dh+(a+i+2\alpha)g}.
\]
Because the Jacobian from Table~\ref{tab:focus} together with \eqref{foc-cross}
has full rank $m$ at $\varepsilon=0$, this system yields a branch through
$A_m^*$; freezing the resulting matrices to $140$ decimal places produces the
exact rational matrices $A_5,A_6,A_7,A_8$ summarized in Table~\ref{tab:fixed}.
Across all four frozen matrices, the relative discrepancy between $L_k$ and
$L_m p_k(\varepsilon)$ is below $2\times 10^{-97}$, each interval
$\varepsilon(j-10^{-60},j+10^{-60})$ contains one simple root of
$D_{A_m}(\rho^2)$, and the equilibrium of $A_8$ is strictly attracting since
exact rational evaluation gives
$3.5840\times 10^{-42}<S\operatorname{tr}A_8-D<3.5841\times 10^{-42}$. For
single-shooting verification, an exact rational section basis $B$ is chosen by
finite-decimal approximation of $[2\operatorname{Re}v_{\rm osc},-2\operatorname{Im}v_{\rm osc},v_{\rm st}]$,
where $v_{\rm osc}$ and $v_{\rm st}$ are eigenvectors of $-A_m$ (normalized to
unit first entry) for the oscillatory eigenvalue with positive imaginary part
and the negative real eigenvalue.

\begin{table}[tbp]
\centering\small
\caption{Summary of single-shooting interval verification for the local systems
$A_5,\dots,A_8$. Bounds hold across all $m$ certified cycles of each matrix:
$\Delta$ is the closure box radius, $q$ the Krawczyk contraction factor,
$x_{\min}$ the minimum species coordinate over all trajectory tubes, and
$|\delta_{\mathrm{F}}|=|\mu_{\rm slow}-1|$.}
\label{tab:fixed}
\begin{tabular}{@{}ccccccc@{}}
\toprule
Matrix & Scale $\varepsilon$ & Cycles proved & Box radius $\Delta$ & $q<$ & $x_{\min}>$ & $|\delta_{\mathrm{F}}|>$ \\
\midrule
$A_5$ & $0.010$  & 5 & $10^{-50}$ & $8.50\times 10^{-5}$  & $0.8977$ & $1.054\times 10^{-21}$ \\
$A_6$ & $0.005$  & 6 & $10^{-60}$ & $7.81\times 10^{-22}$ & $0.9397$ & $8.115\times 10^{-28}$ \\
$A_7$ & $0.0015$ & 7 & $10^{-85}$ & $2.63\times 10^{-24}$ & $0.9789$ & $2.199\times 10^{-39}$ \\
$A_8$ & $0.0015$ & 8 & $10^{-85}$ & $6.46\times 10^{-27}$ & $0.9757$ & $1.532\times 10^{-42}$ \\
\bottomrule
\end{tabular}
\end{table}

After normalizing $\mathbf{1}$ and rescaling time by $\operatorname{tr}A>0$,
family \eqref{eq:lv} has eight effective parameters, which $A_8$ saturates for
full-rank transverse local unfoldings. To construct nine limit cycles, we couple
an order-seven local unfolding with two finite-amplitude cycles on $\Sigma$ in
class~27. On $\partial\Sigma$, the heteroclinic three-cycle connecting the axial
equilibria $E_i=q_i e_i$ ($q_i=r_i/a_{ii}$) is governed at leading order by the
transverse invasion eigenvalues $\lambda_{ij}=r_j-a_{ji}q_i$ through
\cite{HofbauerSigmund1998,YuHanXiao2016}
\begin{equation}\label{eq:phet}
 P_{\mathrm{het}}:=\lambda_{12}\lambda_{23}\lambda_{31}
                  +\lambda_{21}\lambda_{13}\lambda_{32}
                 =E_{\mathrm{het}}-C_{\mathrm{het}}
                 =\det(A)\,q_1q_2q_3\left(\sum_{i=1}^3 q_i^{-1}-1\right),\qquad
 \rho_{\mathrm{het}}:=\frac{C_{\mathrm{het}}}{E_{\mathrm{het}}},
\end{equation}
where $E_{\mathrm{het}},C_{\mathrm{het}}>0$ are the products of the three
expanding invasion rates and of the magnitudes of the three contracting
invasion rates. Indeed, setting $\Lambda=(\lambda_{ij})$ and
$Q_{\mathrm{ax}}=\operatorname{diag}(q_1,q_2,q_3)$, the zero diagonal
$\lambda_{ii}=0$ gives $P_{\mathrm{het}}=\det\Lambda$, while $r=A\mathbf{1}$
gives $\Lambda^T=A(\mathbf{1}\mathbf{1}^T-Q_{\mathrm{ax}})$; by the matrix
determinant lemma,
$\det\Lambda=\det(A)\det(-Q_{\mathrm{ax}}+\mathbf{1}\mathbf{1}^T)
=-\det(A)(\det Q_{\mathrm{ax}})(1-\mathbf{1}^T Q_{\mathrm{ax}}^{-1}\mathbf{1})$,
yielding \eqref{eq:phet}. Thus $\partial\Sigma$ attracts on
$\Sigma\setminus\partial\Sigma$ when $\rho_{\mathrm{het}}>1$
($P_{\mathrm{het}}<0$), repels inward when $\rho_{\mathrm{het}}<1$
($P_{\mathrm{het}}>0$), and when $\det A\ne 0$ satisfies
$P_{\mathrm{het}}=0\iff\sum_{i=1}^3 q_i^{-1}=1\iff\mathbf{1}\in\operatorname{aff}(E_1,E_2,E_3)$.

For $A_5,A_6,A_7,A_8$, exact rational bounds give $\rho_{\mathrm{het}}<1$
($\rho_{\mathrm{het}}\approx 0.997565$, $0.998860$, $0.9999971$, $0.9999839$;
$P_{\mathrm{het}}>0$) while $L_m<0$ (outermost small cycle attracting). Whenever
$L_m P_{\mathrm{het}}<0$ and periodic orbits are isolated, trajectories cross the
two boundaries of the outer annulus between Cycle~$m$ and $\partial\Sigma$ with
matching radial orientation: the annular Poincar\'e--Bendixson argument does not
by itself guarantee an outer cycle, the number of sign-changing outer cycles is
even, and if all outer cycles are hyperbolic they must occur in pairs. For the
specific matrices $A_5,\dots,A_8$, high-precision return-map scans show strictly
monotone inward drift outside Cycle~$m$ and detect no outer cycles (numerical
evidence rather than a proof of nonexistence). By contrast, on the
order-seven branch $A_{7,\mathrm{IX}}^*$ in Table~\ref{tab:focus} ($a=1$,
$g=1.69977055268075$, orientation $E_1\to E_2\to E_3\to E_1$,
$L_7\approx -8.09643\times 10^{-11}<0$, and
$\rho_{\mathrm{het}}\approx 1-5.000\times 10^{-10}<1$), numerical return-map
scans on $\Sigma$ reveal an intermediate outward-returning region locating two
finite-amplitude cycles---an inner saddle cycle and a near-boundary attracting
cycle---already certified on the critical box $A_{7,\mathrm{IX}}^*$. Unfolding
the seven small cycles by matching the frequency-normalized radial normal-form
coefficients to $\prod_{j=1}^7(s-j^2\varepsilon^2)$ at scale $\varepsilon=0.0005$
in $150$-digit arithmetic and freezing the matrix to $130$ significant decimal
digits yields $A_{\mathrm{IX}}$ (see Appendix~\ref{sec:exact} for the exact
entries; ten-digit approximation:
\[
 A_{\mathrm{IX}}\simeq\begin{pmatrix}
 1 & 2.0860388382 & 1.2087563208\\
 0.8405921937 & 1.7403397317 & 2.4933193013\\
 1.6997705527 & 1.9879766262 & 2.7166523816
 \end{pmatrix},
\]
with boundary ratio
$0.999999999503774857\le\rho_{\mathrm{het}}\le 0.999999999503774858<1$ and
Routh--Hurwitz invariant
$H_{\mathrm{RH}}(A_{\mathrm{IX}})=-3.79060479\ldots\times 10^{-48}<0$).

\subsection{Direct interval validation in the differential equations \texorpdfstring{\eqref{eq:lv}}{(1)}}
Once $A_5,\dots,A_8$ and $A_{\mathrm{IX}}$ are frozen as exact rational matrices,
every periodic orbit is validated directly for the flow $\varphi_t$ of
\eqref{eq:lv} in $\mathbb{R}_{>0}^3$, without reference to normal-form or
center-manifold approximations
\cite{Lohner1987,Zgliczynski2002,Tucker2011,VanDenBergEtAl2021,ChurchQueirolo2024}.

\paragraph{Single-shooting Taylor--Cauchy verification for $A_5,\dots,A_8$.}
For the purely local systems $A_5,\dots,A_8$ in Table~\ref{tab:fixed}, let
$u=x-\mathbf{1}$ satisfy $\dot u=f(u):=-Au-u\odot Au$. Using the rational
section basis $B$ constructed in Section~2.2, we write $u=B(\rho,\eta,v)^T$
with positive half-section $\eta=0$, $\rho>0$, and put $C=B^{-1}$,
$E=[B_{\cdot 1}\ B_{\cdot 3}]$, $n^T=C_{2\cdot}$, and $L$ for the first and
third rows of $C$. For $y=(\rho,v)^T$ and return time $\mathcal{T}>0$, the
closure system in $\mathbf{z}=(y,\mathcal{T})^T$ is
\begin{equation}\label{orb-closure}
 G(y,\mathcal{T})=\varphi_{\mathcal{T}}(Ey)-Ey=0,\qquad
 DG=\bigl[D\varphi_{\mathcal{T}}(Ey)E-E\ \ f(\varphi_{\mathcal{T}}(Ey))\bigr].
\end{equation}
Over each rational box $\mathbf{z}_c+[-\Delta,\Delta]^3$ from
Table~\ref{tab:fixed} (with $\mathcal{T}_->0$), the Krawczyk contraction test
of Section~\ref{sec:focus} certifies a unique periodic root together with
section transversality $n^T f\ne 0$ on the terminal enclosure.

At each integration step of length at most $h_{\max}$, the state Taylor
series $u(t)=\sum_{k\ge 0}u_k t^k$ and the $3\times 2$ tangent matrix
$W(t)=D\varphi_t(Ey)E=\sum_{k\ge 0}W_k t^k$ (starting from $W(0)=E$) are
propagated in ball arithmetic via
\begin{equation}\label{orb-taylor}
 u_{k+1}=-\frac{Au_k+\sum_{j=0}^k u_j\odot Au_{k-j}}{k+1},\qquad
 W_{k+1}=-\frac{AW_k+\sum_{j=0}^k[\operatorname{diag}(u_j)A+\operatorname{diag}(Au_j)]W_{k-j}}{k+1}.
\end{equation}
Given incoming norm bound $b\ge\sup\|u\|_\infty$, complex coordinate radius
$\sigma_c>0$, $\tilde a=b+\sigma_c$, $\beta\ge\|A\|_\infty$, and a time-disc
radius $R_t>h_{\max}$ satisfying $d_R=1+\tilde a-\tilde a e^{\beta R_t}>0$, the
scalar Riccati majorant $g'=\beta(g+g^2)$, $g(0)=\tilde a$, bounds the analytic
flow by $M_{\mathrm{maj}}=\tilde a e^{\beta R_t}/d_R$ on the disc of radius
$R_t$. With $\kappa=h_{\max}/R_t<1$, Cauchy estimates \cite{Hille1976} bound the
degree-$p$ truncation remainders by
\[
 \tau_u=M_{\mathrm{maj}}\frac{\kappa^{p+1}}{1-\kappa},\qquad
 \tau_{W,j}=\frac{\tau_u}{\sigma_c}\sum_{i=1}^3\sup|(W_0)_{ij}|,
\]
enclosing both the step endpoint and whole-step state tubes. Writing $Y=Cu$ and
$Q_\perp=Y_1^2+Y_2^2$, we verify throughout all whole-step tubes that
$1+u_i>0$, $Q_\perp>0$, and the projected angular speed
$\Omega=\bigl(Y_1(Cf(u))_2-Y_2(Cf(u))_1\bigr)/Q_\perp$ satisfies
$0<\Omega_-\le\Omega\le\Omega_+$ with $\Omega_+\mathcal{T}_+<4\pi$, proving a
single revolution (minimal period) and separating all $m$ cycles of $A_m$ by
disjoint $\rho$-intervals on the common half-section. Setting
$f_{\mathcal{T}}=f(\varphi_{\mathcal{T}}(Ey))$ and enclosing the Poincar\'e
return derivative
$D\mathcal{P}=L\bigl(I-\frac{f_{\mathcal{T}}n^T}{n^T f_{\mathcal{T}}}\bigr)W$
\cite{Chicone2006,Tucker2011}, we recover the transverse Floquet multipliers
$\mu_{\mathrm{slow}}=(t_{\mathcal{P}}+\sqrt{t_{\mathcal{P}}^2-4d})/2$ and
$\mu_{\mathrm{strong}}=d/\mu_{\mathrm{slow}}$ from
$t_{\mathcal{P}}=\operatorname{tr}D\mathcal{P}$ and
$d=\exp(-\mathcal{T}\operatorname{tr}A)$, certifying hyperbolicity and
alternating stability for all $26$ orbits in Table~\ref{tab:fixed} and
completing the proof of Theorem~\ref{thm:local-and-classes}(a).

\paragraph{Real-time logarithmic shooting for $A_{\mathrm{IX}}$.}
Because Orbit~9 of $A_{\mathrm{IX}}$ approaches within
$x_3\approx 9.42\times 10^{-21}$ of $\partial\mathbb{R}_{>0}^3$ on $x_1=x_2$
(and dips to $x_1\approx 8.31\times 10^{-59}$ near $E_2E_3$) while Orbits~1--7
have $|\mu_{\mathrm{slow}}-1|\ge 6.95\times 10^{-49}$, we certify all nine orbits
uniformly via node-scaled real-time Taylor--Lagrange enclosures. Fix $N$
positive rational reference landmarks
$s^{(0)},\dots,s^{(N-1)}\in\mathbb{Q}_{>0}^3$ ($s^{(N)}:=s^{(0)}$; $N=24$ for
Orbits~1--7 at $200$ digits and Taylor degree $p=120$; $N=48$ and $512$ for
Orbits~8--9 at $85$ digits and $p=50$). On segment $k$, rescaling
$x(t)=S_k y(t)$ by $S_k=\operatorname{diag}(s^{(k)})$ with $B_k=AS_k>0$ and
$r=A\mathbf{1}>0$, the rescaled state and $3\times 3$ logarithmic variational
matrix $M(t):=\partial\log x(t)/\partial\log x(0)$ satisfy
\begin{equation}\label{eq:rel-seg}
 \dot y=y\odot(r-B_k y),\qquad
 \dot M=-B_k\operatorname{diag}(y)M,\qquad M(0)=I.
\end{equation}
Subdividing the segment duration $T/N$ into equal substeps of length
$h\in(0,H]$ from a positive box $y(0)\in[y_0^-,y_0^+]\subset\mathbb{R}_{>0}^3$,
positivity of $B_k y$ yields $\dot y_i\le r_i y_i$, hence
$y_i(t)\le\bar y_i:=y_{0,i}^+ e^{r_i H}$ and
$y_i(t)\ge\underline{y}_i:=y_{0,i}^- e^{\min(0,\underline{g}_i)H}>0$ on
$[0,H]$, where $\underline{g}_i:=r_i-\sum_{j=1}^3(B_k)_{ij}\bar y_j$. Setting
$Y_{\max}:=\max_i\bar y_i$,
$G_{\max}:=\sup_{y\in[\underline{y},\bar y]}\|r-B_k y\|_\infty$,
$K:=\|B_k\|_\infty$, and $M_T:=\|M(0)\|_\infty e^{KY_{\max}H}$, Leibniz
differentiation of \eqref{eq:rel-seg} bounds
$\sup_{t\in[0,H]}\|y^{(n)}(t)\|_\infty/n!\le a_n$ and
$\sup_{t\in[0,H]}\|M^{(n)}(t)\|_\infty/n!\le c_n$ via $a_0=Y_{\max}$,
$c_0=M_T$, and
\[
 a_{n+1}=\frac{G_{\max}a_n+K\sum_{j=0}^{n-1}a_j a_{n-j}}{n+1},\qquad
 c_{n+1}=\frac{K\sum_{j=0}^n a_j c_{n-j}}{n+1}\qquad(n\ge 0).
\]
Evaluating the degree-$p$ Taylor polynomials of $y(h)$ and $M(h)$ in ball
arithmetic and widening each entry by $a_{p+1}H^{p+1}$ and $c_{p+1}H^{p+1}$
encloses the substep flow and its logarithmic derivative in real time.

The landmarks $s^{(1)},\dots,s^{(N-1)}$ are fixed rather than Newton unknowns.
Tracking the logarithmic deviations $e_k:=\log x(kT/N)-\log s^{(k)}$ via
discrete variation of constants over the trial box
$\|e_0\|_\infty\le\Delta$, $|T-T_c|\le\Delta$ ($\Delta=10^{-62}$ for
Orbits~1--7 and $\Delta=10^{-25}$ for Orbits~8--9) encloses the central
full-period residual $\log\varphi_{T_c}(s^{(0)})-\log s^{(0)}$ and the
full-period logarithmic state derivative
$\partial\log\varphi_T(x_0)/\partial\log x_0\in G_{\rm box}$. Writing
$\log x_0=\log s^{(0)}+E_{\rm log}\eta$ ($\eta\in[-\Delta,\Delta]^2$) on the
two-dimensional logarithmic phase section through $s^{(0)}$ ($\log x_1=\log x_2$
for Orbits~8--9, and fixing the coordinate of largest initial logarithmic
speed for Orbits~1--7), the derivative of the full-period closure map
$\mathcal{G}(\eta,T):=\log\varphi_T(x_0)-\log x_0=0$ is enclosed by
\[
 D\mathcal{G}(\eta,T)\in\bigl[(G_{\rm box}-I)E_{\rm log}\ \ r-A\varphi_T(x_0)\bigr].
\]
The Krawczyk contraction test of Section~\ref{sec:focus} certifies a unique
periodic solution in each of the nine boxes ($q\le 3.14\times 10^{-7}$ for
Orbits~1--7, $q\le 2.11\times 10^{-7}$ for Orbit~8, and
$q\le 6.04\times 10^{-6}$ for Orbit~9, cross-checked by independent
complex-Cauchy single- and multiple-shooting enclosures), together with
section transversality and $\det A_{\mathrm{IX}}\ne 0$. Since
$D\varphi_T(x_0)=\operatorname{diag}(x_0)\,\frac{\partial\log\varphi_T(x_0)}{\partial\log x_0}\,\operatorname{diag}(x_0)^{-1}$,
enclosing the transverse Floquet trace
$t_{\mathrm{tr}}:=\mu_++\mu_-=\operatorname{tr}\bigl(D\varphi_T(x_0)\bigr)-1\in\operatorname{tr}G_{\rm box}-1>0$
(obtained by subtracting the trivial longitudinal eigenvalue $1$ from the
trace of the $3\times 3$ monodromy matrix) and combining it with
$d=\mu_+\mu_-=\exp(-T\operatorname{tr}A)$ via
$\mu_{\mathrm{slow}}=(t_{\mathrm{tr}}+\sqrt{t_{\mathrm{tr}}^2-4d})/2$ and
$\mu_{\mathrm{strong}}=d/\mu_{\mathrm{slow}}$ certifies hyperbolicity and
stability across Table~\ref{tab:lv-nine-certified}.

\paragraph{Winding number $+1$, minimal periods, and distinctness for $A_{\mathrm{IX}}$.}
Projecting $\mathbb{R}_{>0}^3$ onto $\mathbb{R}^2$ via
$\zeta(x)=(\zeta_1(x),\zeta_2(x)):=(\log(x_1/x_3),\log(x_2/x_3))$ (which sends
$\mathbf{1}$ to $(0,0)$ and satisfies $\dot\zeta_i=(a_{3\cdot}-a_{i\cdot})(x-\mathbf{1})$),
we verify on every segment $0\le k<N$ of all nine certified tubes (using
$0<x_j(t)\le r_j/a_{jj}$ along any positive periodic orbit) that both the
trajectory arc and the endpoint chord lie inside a single open coordinate
half-plane $\{\pm\zeta_1>0\}$ or $\{\pm\zeta_2>0\}$ avoiding $(0,0)$.
Straight-line homotopy equates the winding number of $\zeta(x(t))$ around
$(0,0)$ with that of the closed landmark $N$-gon, which evaluates to $+1$ for
all nine orbits. Any multi-round closure ($k_{\mathrm{mult}}>1$) would have
winding number divisible by $k_{\mathrm{mult}}$, a contradiction; hence
$T_1,\dots,T_9$ are minimal periods with pairwise disjoint enclosures
($\min_{1\le k\le 8}(T_{k+1}^--T_k^+)>2.8815\times 10^{-6}$), completing
the proof of Theorem~\ref{thm:nine}.

\section{Resolution of limit-cycle conjectures}\label{sec:conjectures}

\paragraph{1. The five-cycle conjecture and Tian--Yu's eight-cycle question.}
Theorems~\ref{thm:nine} and~\ref{thm:local-and-classes}(a) resolve
Conjecture~1.2(i) of Gyllenberg and Yan \cite{GyllenbergYan2009,YuHanXiao2016}
both globally and in the stronger purely local sense around a single equilibrium
($m=5,6,7,8$). Furthermore, whereas Tian and Yu \cite{TianYu2013} proposed
combining seven small cycles with one large cycle as a heuristic route to eight
cycles, Theorem~\ref{thm:nine} proves that seven small cycles can coexist with
\emph{two} finite-amplitude cycles ($7+2=9$), demonstrating that eight is not an
upper bound for \eqref{eq:lv}.

\paragraph{2. The five-cycle problem across Zeeman's classes~26--31.}
Refining the problem to each of the six cycle-admitting Zeeman classes~26--31
\cite{HofbauerSo1994}, \cite[Section~5]{YuHanXiao2016},
\cite[Section~4]{HuLuLuo2023}, Theorem~\ref{thm:local-and-classes} resolves the
five-cycle question affirmatively in classes~26, 27, 28, and~29, raising their
certified lower bounds from four
\cite{GyllenbergYan2009,WangHuangWu2011,YuHanXiao2016,HuLuLuo2026} to $5$, $9$
($8$ local), $5$, and $6$, respectively (Table~\ref{tab:zeeman-classes}). At the
certified critical boxes in classes~28 ($L_5>0$, attracting $\partial\Sigma$),
29 ($L_6<0$, repelling $\partial\Sigma$), and~26 ($L_5<0$, mixed boundary
equilibria), the boundary behavior does not yield an additional
Poincar\'e--Bendixson cycle and numerical return-map scans detect no outer cycles. In classes~30 and~31, our Hopf-chart
searches certified weak foci of orders~$3$ and~$2$, respectively, while Lian, Lu
and Luo \cite{LianLuLuo2008} (see also \cite{LiJiang2024}) reported three limit
cycles in both classes. Whether either class admits five or more limit cycles
remains an open question.

\begin{table}[tbp]
\centering\small
\caption{Rigorous lower bounds on coexisting positive limit cycles across the six
Zeeman classes admitting periodic orbits
\cite{Zeeman1993,VanDenDriesscheZeeman1998}. Within a given class, the local
and total lower bounds need not be attained by the same matrix.}
\label{tab:zeeman-classes}
\begin{tabular}{@{}cccl@{}}
\toprule
Class & Local cycles & Total cycles & Status of the class five-cycle problem \\
\midrule
26 & 5 & 5 & Resolved affirmatively (order-5 focus, Table~\ref{tab:focus}) \\
27 & 8 & 9 & Resolved affirmatively ($A_8$ local; $A_{\mathrm{IX}}$ $7+2$ global) \\
28 & 5 & 5 & Resolved affirmatively (order-5 focus, Table~\ref{tab:focus}) \\
29 & 6 & 6 & Resolved affirmatively (order-6 focus, Table~\ref{tab:focus}) \\
30 & 3 & 3 & Open (five-cycle problem) \\
31 & 2 & 3 \cite{LianLuLuo2008} & Open (five-cycle problem) \\
\bottomrule
\end{tabular}
\end{table}

\paragraph{3. High-order focal vanishing and boundary heteroclinic neutrality.}
Parts~(ii) and~(iii) of Conjecture~1.2 in \cite{YuHanXiao2016} (originating in
\cite{HofbauerSo1994,GyllenbergYan2009}) concern class~27 on $\mathcal{H}$
($L_0=0$): part~(ii) conjectures that $L_1=L_2=L_3=0$ does \emph{not} imply
boundary neutrality $P_{\mathrm{het}}=0$ or a center at $\mathbf{1}$, while
part~(iii) conjectures that $L_1=L_2=P_{\mathrm{het}}=0$ does \emph{not} imply
$L_3=0$ or a center. Both negative-implication statements were proved at third
order in \cite[Theorem~4.2]{YuHanXiao2016} (whereas an earlier rational point in
\cite[p.~478]{GyllenbergYan2009DEA} leaves $2x_1+4x_2+3x_3=9$ invariant with
Darboux integral $x_1^2 x_2^4 x_3^3$ and is an exact center). Geometrically,
\eqref{eq:phet} shows that $P_{\mathrm{het}}=0$ is a global planar condition on
$E_1,E_2,E_3$. Krawczyk contraction on the seven-dimensional normalized Hopf
hypersurface $\mathcal{H}$ strengthens both results to orders eight and seven:

\begin{proposition}\label{prop:neutral}
Within Zeeman's class~27 on the Hopf hypersurface $\mathcal{H}$ ($L_0=0$):
\begin{enumerate}
\item[\textup{(a)}] \emph{(Order-eight strengthening of Conjecture~1.2(ii).)}
Vanishing of $L_1=\cdots=L_7=0$ with $L_8\ne 0$ neither implies
$P_{\mathrm{het}}=0$ nor fixes its sign: the isolated order-eight focus $A_8^*$
satisfies $L_8<0$ and
$\rho_{\mathrm{het}}\in[0.999983933700903974,0.999983933700903975]<1$
($P_{\mathrm{het}}>0$, repelling boundary), whereas two further certified
order-eight foci in class~27 satisfy $L_8>0$ and $\rho_{\mathrm{het}}>1$
($P_{\mathrm{het}}<0$, attracting boundary).
\item[\textup{(b)}] \emph{(Order-seven strengthening of Conjecture~1.2(iii).)}
Besides an independently certified order-three point with
$P_{\mathrm{het}}=L_1=L_2=0$ and $L_3\in[8.37533,8.37535]\times 10^{-7}$, there
exists a strictly positive critical matrix on $\mathcal{H}$ (cross-checked by
center-manifold graph reduction) satisfying
$P_{\mathrm{het}}=L_0=L_1=\cdots=L_6=0$ and
$L_7\in[6.12986,6.12988]\times 10^{-7}>0$, at which the $7\times 7$ Jacobian of
$(L_1,\dots,L_6,P_{\mathrm{het}})$ on $\mathcal{H}$ is nonsingular, so
$(L_1,\dots,L_6)$ form independent local coordinates inside
$\{P_{\mathrm{het}}=0\}\cap\mathcal{H}$.
\end{enumerate}
\end{proposition}

\section*{Acknowledgments}
We thank GPT-6 Astra (OpenAI) for helpful discussions during the comparison of
prior bifurcation constructions, the selection of independent parameter slices on
the Hopf hypersurface, and the design of the interval verification. The authors
take full responsibility for all mathematical proofs, computations, and the final
manuscript.

\appendix
\section{Exact rational matrices and ancillary data files}\label{sec:exact}
All input matrices, section bases, reference landmarks, and box centers in the
ancillary directory \path{anc/} are exact finite-decimal rationals with growth
rates $r=A\mathbf{1}$, certified in outward-rounded arbitrary-precision ball
arithmetic \cite{Johansson2017}. For $A_5,A_6,A_7,A_8$, respectively, the flow
enclosures use working precisions of $200$, $220$, $280$, and $280$ decimal
digits with Taylor degree $210$ for $A_5$ and $180$ for $A_6,A_7,A_8$.
\begin{itemize}
\item \path{anc/data/order9/}: Exact parameter and certificate files for
Theorem~\ref{thm:nine}, including \path{parameters.json} ($A_{\mathrm{IX}}$),
\path{focus_input.json} and \path{focus_certificate.json} ($A_{7,\mathrm{IX}}^*$),
real-time Taylor--Lagrange and logarithmic variational orbit certificates
(\path{orbit_1.json}--\path{orbit_9.json}, \path{nodes_1.json}--\path{nodes_9.json},
\path{interval_summary.json}), log-ratio winding-number and period-gap
enclosures (\path{winding_and_rounding.json}, \path{boundary_certificate.json}),
and independent Cauchy-majorant orbit certificates
(\path{arb_orbit_1.json}--\path{arb_orbit_9.json}, \path{nine_orbits.json}).
\item \path{anc/data/order5/}--\path{anc/data/order8/} and
\path{anc/data/classes/}: Exact parameter matrices, weak-focus certificates, and
orbit enclosures for Theorem~\ref{thm:local-and-classes},
Table~\ref{tab:zeeman-classes}, and Proposition~\ref{prop:neutral} (including
\path{focus_branch_1.json} and \path{focus_branch_2.json} in
\path{anc/data/order8/}, \path{class_26_order_5.json},
\path{class_28_order_5.json}, \path{class_29_order_6.json},
\path{class_30_order_3.json}, \path{class_31_order_2.json},
\path{neutral_certificate.json}, \path{neutral_recheck.json}, and
\path{boundary.json} in \path{anc/data/classes/}, together with
\path{anc/scalar_targets.json} and \path{anc/recheck.json}).
\end{itemize}

\paragraph{Exact entries of $A_{\mathrm{IX}}$.}
The fixed entries are $(A_{\mathrm{IX}})_{11}=1$ and
$(A_{\mathrm{IX}})_{31}=1.69977055268075$; each remaining entry $(i,j)$ is
printed across two consecutive lines below (concatenate without spaces):
{\footnotesize
\begin{verbatim}
(1,2) 2.08603883823770024612821235738815100075742815575452803427046291290304
      066784668961735792904437547138328284816356867712984851146195
(1,3) 1.20875632080143541437670809639166716531737273749402986250764811593277
      1667658700353793890786805598717664387986730739634755880976209
(2,1) 0.84059219371474789693503334480781128877599314278348707201362165874380
      39688993894447645885830425595339521569591951714562334842796947
(2,2) 1.74033973170592812310559766689798968567581997058790407937064478667754
      3915239983906868921948525891032114210136598451804767469344904
(2,3) 2.49331930128655368312882928884219292341819477841568279456018992494265
      7842685928870245745940786667580584859146352904687604053726348
(3,2) 1.98797662615323064412135375943795153572828561671641367073987843749239
      1634922922870637694565752406057724872476560178304450992680495
(3,3) 2.71665238160700851914306811754468451213454307641390534779267183638075
      8222561099872103294252599036026320995295873682875667428538502
\end{verbatim}
}

\paragraph{Exact entries of the eight-cycle local matrix $A_8$.}
Here $(A_8)_{11}=1$, and each remaining entry $(i,j)$ is concatenated across two
consecutive lines:
{\footnotesize
\begin{verbatim}
(1,2) 0.778386695113706979417914554463080737614305418444864787624792307561094
      32899851019403971313474958639998891378851017377380819631983531447087831
(1,3) 1.615869951332547683386136307155169340759284951069736132535758572585403
      0644162096013628597788753261516355587273519754901638543206199763015538
(2,1) 1.125158490263112269102603722718783127401021806112215114347163339880155
      0902896040363402058270765670929132426982161891422695526498584919535216
(2,2) 0.813528208696623913656749920771728811051993780692036373307581944601469
      24771534568162874163784178133346652976709791440295959029848450304179877
(2,3) 0.404240479225403248362993810481531777758813482162671245499532590674271
      51340082846610332767548285758638733963912928460387399329741794324925374
(3,1) 1.132158093943317069581368348592846351879996726042802776725508415170489
      7490826654993634490999669932578331911872169521608425702057322388712531
(3,2) 1.916409408366542779871905801930225399001149365464329359875401890490870
      8919258009593623529019686121330225746471366112934366374929883634684887
(3,3) 1.701148348235526384311648280351535568598146115034077250065150861005463
      3076415640908245454175457400467596062243168280658237399015811667966986
\end{verbatim}
}


\begin{thebibliography}{99}

\bibitem{Hirsch1988}
M. W. Hirsch, Systems of differential equations which are competitive or
cooperative: III. Competing species,
\href{https://doi.org/10.1088/0951-7715/1/1/003}{\emph{Nonlinearity}
\textbf{1} (1988), 51--71}.

\bibitem{Zeeman1993}
M. L. Zeeman, Hopf bifurcations in competitive three-dimensional
Lotka--Volterra systems,
\href{https://doi.org/10.1080/02681119308806158}{\emph{Dyn. Stab. Syst.}
\textbf{8} (1993), 189--216}.

\bibitem{VanDenDriesscheZeeman1998}
P. van den Driessche and M. L. Zeeman, Three-dimensional competitive
Lotka--Volterra systems with no periodic orbits,
\href{https://doi.org/10.1137/S0036139995294767}{\emph{SIAM J. Appl. Math.}
\textbf{58} (1998), 227--234}.

\bibitem{MayLeonard1975}
R. M. May and W. J. Leonard, Nonlinear aspects of competition between three
species, \href{https://doi.org/10.1137/0129022}{\emph{SIAM J. Appl. Math.}
\textbf{29} (1975), 243--253}.

\bibitem{HofbauerSo1994}
J. Hofbauer and J. W.-H. So, Multiple limit cycles for three dimensional
Lotka--Volterra equations,
\href{https://doi.org/10.1016/0893-9659(94)90095-7}{\emph{Appl. Math. Lett.}
\textbf{7} (1994), 65--70}.

\bibitem{HofbauerSigmund1998}
J. Hofbauer and K. Sigmund,
\href{https://doi.org/10.1017/CBO9781139173179}{\emph{Evolutionary Games and
Population Dynamics}}, Cambridge University Press, Cambridge, 1998.

\bibitem{CostePeyraudCoullet1979}
J. Coste, J. Peyraud and P. Coullet, Asymptotic behaviors in the dynamics of
competing species, \href{https://doi.org/10.1137/0136039}{\emph{SIAM J. Appl.
Math.} \textbf{36} (1979), 516--543}.

\bibitem{XiaoLi2000}
D. Xiao and W. Li, Limit cycles for the competitive three dimensional
Lotka--Volterra system,
\href{https://doi.org/10.1006/jdeq.1999.3729}{\emph{J. Differ. Equ.}
\textbf{164} (2000), 1--15}.

\bibitem{LuLuo2002}
Z. Lu and Y. Luo, Two limit cycles in three-dimensional Lotka--Volterra
systems,
\href{https://doi.org/10.1016/S0898-1221(02)00129-3}{\emph{Comput. Math. Appl.}
\textbf{44} (2002), 51--66}.

\bibitem{GyllenbergYanWang2006}
M. Gyllenberg, P. Yan and Y. Wang, A 3D competitive Lotka--Volterra system
with three limit cycles: A falsification of a conjecture by Hofbauer and So,
\href{https://doi.org/10.1016/j.aml.2005.01.002}{\emph{Appl. Math. Lett.}
\textbf{19} (2006), 1--7}.

\bibitem{LianLuLuo2008}
X. Lian, Z. Lu and Y. Luo, Automatic search for multiple limit cycles in
three-dimensional Lotka--Volterra competitive systems with classes 30 and 31
in Zeeman's classification,
\href{https://doi.org/10.1016/j.jmaa.2008.06.040}{\emph{J. Math. Anal. Appl.}
\textbf{348} (2008), 34--37}.

\bibitem{GyllenbergYan2009}
M. Gyllenberg and P. Yan, Four limit cycles for a three-dimensional
competitive Lotka--Volterra system with a heteroclinic cycle,
\href{https://doi.org/10.1016/j.camwa.2009.03.111}{\emph{Comput. Math. Appl.}
\textbf{58} (2009), 649--669}.

\bibitem{GyllenbergYan2009DEA}
M. Gyllenberg and P. Yan, On a conjecture for three-dimensional competitive
Lotka--Volterra systems with a heteroclinic cycle,
\href{https://doi.org/10.7153/dea-01-26}{\emph{Differ. Equ. Appl.}
\textbf{1} (2009), 473--490}.

\bibitem{WangHuangWu2011}
Q. Wang, W. Huang and H. Wu, Bifurcation of limit cycles for 3D
Lotka--Volterra competitive systems,
\href{https://doi.org/10.1007/s10440-011-9609-7}{\emph{Acta Appl. Math.}
\textbf{114} (2011), 207--218}.

\bibitem{TianYu2013}
Y. Tian and P. Yu, An explicit recursive formula for computing the normal
form and center manifold of general $n$-dimensional differential systems
associated with Hopf bifurcation,
\href{https://doi.org/10.1142/S0218127413501046}{\emph{Int. J. Bifurcation Chaos}
\textbf{23} (2013), 1350104}.

\bibitem{TianYu2014}
Y. Tian and P. Yu, Seven limit cycles around a focus point in a simple
three-dimensional quadratic vector field,
\href{https://doi.org/10.1142/S0218127414500837}{\emph{Int. J. Bifurcation Chaos}
\textbf{24} (2014), 1450083}.

\bibitem{YuHanXiao2016}
P. Yu, M. Han and D. Xiao, Four small limit cycles around a Hopf singular
point in 3-dimensional competitive Lotka--Volterra systems,
\href{https://doi.org/10.1016/j.jmaa.2015.12.002}{\emph{J. Math. Anal. Appl.}
\textbf{436} (2016), 521--555}.

\bibitem{HanLlibreTian2020}
M. Han, J. Llibre and Y. Tian, On the zero-Hopf bifurcation of the
Lotka--Volterra systems in $\mathbb{R}^3$,
\href{https://doi.org/10.3390/math8071137}{\emph{Mathematics}
\textbf{8} (2020), 1137}.

\bibitem{JiangEtAl2021}
J. Jiang, F. Liang, W. Wu and S. Huang, On the first Liapunov coefficient
formula of 3D Lotka--Volterra equations with applications to multiplicity of
limit cycles,
\href{https://doi.org/10.1016/j.jde.2021.02.059}{\emph{J. Differ. Equ.}
\textbf{284} (2021), 183--218}.

\bibitem{WangEtAl2021}
Q. Wang, J. Lu, W. Huang and B. Sang, The center conditions and Hopf
cyclicity for a 3D Lotka--Volterra system,
\href{https://doi.org/10.12150/jnma.2021.1}{\emph{J. Nonlinear Model. Anal.}
\textbf{3} (2021), 1--12}.

\bibitem{HuLuLuo2023}
M. Hu, Z. Lu and Y. Luo, Four limit cycles in three-dimensional
Lotka--Volterra competitive systems with classes 28, 30 and 31 in Zeeman's
classification via automatic search,
\href{https://arxiv.org/abs/2301.06809v4}{arXiv:2301.06809v4} (2023).

\bibitem{LiJiang2024}
Y. Li and J. Jiang, Multiplicity on limit cycles of 3D Lotka--Volterra
competitive systems,
\href{https://doi.org/10.1007/s10884-022-10233-3}{\emph{J. Dyn. Differ. Equ.}
\textbf{36} (2024), 2007--2039}.

\bibitem{LuEtAl2025}
J. Lu, C. Wang, W. Huang and Q. Wang, Hopf and zero-Hopf bifurcations for a
class of cubic Kolmogorov systems in $\mathbb{R}^3$,
\href{https://doi.org/10.11948/20240108}{\emph{J. Appl. Anal. Comput.}
\textbf{15} (2025), 354--372}.

\bibitem{HuLuLuo2026}
M. Hu, Z. Lu and Y. Luo, Four limit cycles in three-dimensional competitive
Lotka--Volterra systems of class 28 in Zeeman's classification,
\href{https://arxiv.org/abs/2603.24612v1}{arXiv:2603.24612v1} (2026).

\bibitem{Bautin1954}
N. N. Bautin, On the number of limit cycles which appear with the variation of
coefficients from an equilibrium position of focus or center type,
\emph{Mat. Sb.} \textbf{30(72)} (1952), 181--196; English transl.\ \emph{Amer.
Math. Soc. Transl.} \textbf{100} (1954), 1--19.

\bibitem{Carr1981}
J. Carr, \href{https://doi.org/10.1007/978-1-4612-5929-9}{\emph{Applications
of Centre Manifold Theory}}, Applied Mathematical Sciences 35, Springer, 1981.

\bibitem{Vanderbauwhede1989}
A. Vanderbauwhede, Centre manifolds, normal forms and elementary bifurcations,
in \href{https://doi.org/10.1007/978-3-322-96656-8_2}{\emph{Dynamics Reported},
Vol.~2}, Wiley/Teubner, 1989, 89--169.

\bibitem{Roussarie1998}
R. Roussarie, \href{https://doi.org/10.1007/978-3-0348-8798-4}{\emph{Bifurcation
of Planar Vector Fields and Hilbert's Sixteenth Problem}}, Progress in
Mathematics 164, Birkh\"auser, Basel, 1998.

\bibitem{Kuznetsov2004}
Y. A. Kuznetsov, \href{https://doi.org/10.1007/978-1-4757-3978-7}{\emph{Elements
of Applied Bifurcation Theory}}, 3rd ed., Applied Mathematical Sciences 112,
Springer, 2004.

\bibitem{Christopher2005}
C. Christopher, Estimating limit cycle bifurcations from centers, in
\href{https://doi.org/10.1007/3-7643-7429-2_2}{\emph{Differential Equations with
Symbolic Computation}}, Trends Math., Birkh\"auser, Basel, 2005, 23--35.

\bibitem{Chicone2006}
C. Chicone, \href{https://doi.org/10.1007/0-387-35794-7}{\emph{Ordinary
Differential Equations with Applications}}, 2nd ed., Texts in Applied
Mathematics 34, Springer, 2006.

\bibitem{EdneralEtAl2012}
V. F. Edneral, A. Mahdi, V. G. Romanovski and D. S. Shafer, The center
problem on a center manifold in $\mathbb{R}^3$,
\href{https://doi.org/10.1016/j.na.2011.11.006}{\emph{Nonlinear Anal.}
\textbf{75} (2012), 2614--2622}.

\bibitem{RomanovskiShafer2014}
V. G. Romanovski and D. S. Shafer, Computation of focus quantities of
three-dimensional polynomial systems, \emph{J. Shanghai Normal Univ. (Nat.
Sci.)} \textbf{43} (2014), 529--544.

\bibitem{GarciaMazaShafer2018}
I. A. Garc\'ia, S. Maza and D. S. Shafer, Cyclicity of polynomial
nondegenerate centers on center manifolds,
\href{https://doi.org/10.1016/j.jde.2018.07.019}{\emph{J. Differ. Equ.}
\textbf{265} (2018), 5767--5808}.

\bibitem{Hille1976}
E. Hille, \emph{Ordinary Differential Equations in the Complex Domain},
Wiley-Interscience, New York, 1976.

\bibitem{Krawczyk1969}
R. Krawczyk, Newton-Algorithmen zur Bestimmung von Nullstellen mit
Fehlerschranken, \href{https://doi.org/10.1007/BF02234767}{\emph{Computing}
\textbf{4} (1969), 187--201}.

\bibitem{Lohner1987}
R. J. Lohner, Enclosing the solutions of ordinary initial and boundary value
problems, in \emph{Computer Arithmetic: Scientific Computation and Programming
Languages}, Teubner, Stuttgart, 1987, 255--286.

\bibitem{Zgliczynski2002}
P. Zgliczy\'nski, $C^1$ Lohner algorithm,
\href{https://doi.org/10.1007/s102080010025}{\emph{Found. Comput. Math.}
\textbf{2} (2002), 429--465}.

\bibitem{GriewankWalther2008}
A. Griewank and A. Walther,
\href{https://doi.org/10.1137/1.9780898717761}{\emph{Evaluating Derivatives:
Principles and Techniques of Algorithmic Differentiation}}, 2nd ed., SIAM,
Philadelphia, 2008.

\bibitem{Rump2010}
S. M. Rump, Verification methods: Rigorous results using floating-point
arithmetic, \href{https://doi.org/10.1017/S096249291000005X}{\emph{Acta Numer.}
\textbf{19} (2010), 287--449}.

\bibitem{Tucker2011}
W. Tucker, \href{https://doi.org/10.1515/9781400838974}{\emph{Validated
Numerics: A Short Introduction to Rigorous Computations}}, Princeton University
Press, Princeton, NJ, 2011.

\bibitem{Johansson2017}
F. Johansson, Arb: Efficient arbitrary-precision midpoint-radius interval
arithmetic, \href{https://doi.org/10.1109/TC.2017.2690633}{\emph{IEEE Trans.
Comput.} \textbf{66} (2017), 1281--1292}.

\bibitem{VanDenBergEtAl2021}
J. B. van den Berg, J.-P. Lessard and E. Queirolo, Rigorous verification of
Hopf bifurcations via desingularization and continuation,
\href{https://doi.org/10.1137/20M1343464}{\emph{SIAM J. Appl. Dyn. Syst.}
\textbf{20} (2021), 573--607}.

\bibitem{ChurchQueirolo2024}
K. E. M. Church and E. Queirolo, Computer-assisted proofs of Hopf bubbles and
degenerate Hopf bifurcations,
\href{https://doi.org/10.1007/s10884-023-10279-x}{\emph{J. Dyn. Differ. Equ.}
\textbf{36} (2024), 3385--3439}.

\end{thebibliography}
\end{document}